\documentclass[12pt]{amsproc}
\usepackage[cp1251]{inputenc}
\usepackage[T2A]{fontenc} 
\usepackage[russian]{babel} 
\usepackage{amsbsy,amsmath,amssymb,amsthm,graphicx,color} 
\usepackage{amscd}
\def\le{\leqslant}
\def\ge{\geqslant}

\newtheorem{prop}{Предложение}

\theoremstyle{definition}

\theoremstyle{remark}

\begin {document}
\centerline{УДК 515.162.8+515.164.15}
\unitlength=1mm
\title[Ряды Пуанкаре групп Клейна]
{Ряды Пуанкаре групп Клейна, многочлены Кокстера, представление Бурау и инварианты Милнора}
\author{Г. Г. Ильюта}
\email{ilyuta@mccme.ru}
\address{Московский Государственный Гуманитарный Университет}
\thanks{Работа поддержана грантами РФФИ-07-01-00593, INTAS-05-7805  и НШ-4719.2006.1.}
\maketitle
\begin{abstract} 

 Мы получим несколько формул для определённых Б.~Костантом рядов Пуанкаре групп Клейна (бинарных полиэдральных групп) и многочленов Кокстера (характеристических многочленов монодромии в случае особенностей). Некоторые из них -- обобщение формулы Эбелинга, тождество Кристоффеля-Дарбу и комбинаторная формула -- являются следствиями известных утверждений о характеристическом многочлене графа. Отношения рядов Пуанкаре и многочленов Кокстера представлены в виде ветвящихся цепных дробей, которые являются $q$-аналогами цепных дробей, появляющихся в теории разрешений особенностей и исчислении Кирби (возникает вопрос о существовании $q$-аналогов этих теорий -- возможно, такие обобщения могли бы основываться на результатах А.~Б.~Гивенталя о $q$-монодромии). Оставшиеся формулы связывают отношения некоторых рядов Пуанкаре и многочленов Кокстера с представлением Бурау и инвариантами Милнора стринг-зацеплений -- обобщённых кос, в которых нитям разрешается не быть монотонными. Результаты С.~М.~Гусейн-Заде, Ф.~Дельгадо и А.~Кампильо позволяют рассматривать эти факты как утверждения о рядах Пуанкаре колец функций на особенностях кривых. Использованные в доказательствах результаты о факторизации многочлена Александера-Конвея зацепления и формулы этой статьи подсказывают гипотезу: отношение рядов Пуанкаре колец функций для близких (по примыканию или по расположению в серии) особенностей кривых определяется представлением Бурау или инвариантами Милнора стринг-зацепления, которое является промежуточным объектом при перестройке узла одной особенности в узел другой.

\end{abstract}
 
\section {Введение}
  
  В.~Эбелинг нашёл формулу, представляющую ряд Пуанкаре для инвариантов группы Клейна в виде отношения многочленов Кокстера евклидовой и аффинной диаграмм Дынкина, отвечающих группе согласно соответствию Маккея (исключение составляет аффинная диаграмма Дынкина $\tilde A_{2m}$) \cite{7}. Ранее были известны результаты С.~М.~Гусейн-Заде, Ф.~Дельгадо и А.~Кампильо, связывающие характеристический многочлен монодромии (многочлен Кокстера) с рядом Пуанкаре кольца функций на особенности \cite{3}. Наш подход к теме связей рядов Пуанкаре и многочленов Кокстера основан на аналогии между алгоритмом Евклида и рекуррентными соотношениями для рядов Пуанкаре и многочленов Кокстера или, в более общем контексте, соотношениями для блочных определителей, известными как формула дополнения по Шуру. Другой подход опирается на результаты о факторизации многочлена Александера-Конвея зацепления \cite{14}, \cite{19}.

 Доказательства многих связанных с поверхностью Зейферта, матрицей Александера и многочленом Александера-Конвея результатов теории узлов сводятся к применению формулы дополнения по Шуру (часто формула явно не используется, но её применение может упростить доказательство). В простейшем варианте эта формула позволяет вычислять определитель блочной $2\times 2$-матрицы
$$
 \det {(A_{ij})}=\det {(A_{11}-A_{12}A_{22}^{-1}A_{21})}\det {A_{22}}.   
$$
В теории циклических накрытий дополнения к узлу и её обобщениях формула появляется как тождество для блочного обобщения циркулянта или стандартного матричного представления циклической перестановки $(i\to i+1)\mod n$  \cite{22}, \cite{25}, \cite{26}. В работах, связанных с факторизацией многочлена Александера-Конвея, формула дополнения по Шуру используется в ситуации, когда матрица формы Зейферта одного зацепления определяется частью матрицы формы Зейферта другого зацепления \cite{14}, \cite{19}, \cite{30}. В \cite{4} с помощью этой формулы получены соотношения для многочлена Александера-Конвея и характеристических многочленов операторов из группы монодромии особенности. Мы рассматриваем формулу дополнения по Шуру как первое "деление с остатком"  для некоторого обобщённого алгоритма Евклида. Особенно близкие к алгоритму Евклида формулы получаются, если в качестве $A_{11}$ взять подматрицу порядка $1$.  Примененяя формулу дополнения по Шуру к трёхдиагональной матрице Якоби, можно доказать известные рекуррентные соотношения для ортогональных многочленов (в этом случае достаточно разложений определителя по строкам и столбцам). Аналогия между алгоритмом Евклида и теорией ортогональных многочленов использовалась в \cite{18}, где, в частности, доказан аналог тождества Кристоффеля-Дарбу для безутиана и вронскиана пары многочленов. В \cite{9} имеется соответствующее тождество для характеристического многочлена графа. В п.~2 мы докажем тождество Кристоффеля-Дарбу для многочленов Кокстера, а также аналоги некоторых других формул для алгоритма Евклида, в частности, получим разложения в ветвящиеся цепные дроби. Для евклидовой диаграммы Дынкина $A_n$ цепная дробь равна частному двух соседних многочленов Чебышёва. Для аффинных диаграмм Дынкина $\tilde A_n$, $\tilde D_n$, $\tilde E_6$, $\tilde E_7$, $\tilde E_8$ ветвящиеся цепные дроби выписаны явно и в п.~3 показано, что они равны умноженным на $q$ рядам Пуанкаре для инвариантов групп Клейна. В теории ортогональных многочленов цепные дроби (а также тождество Кристоффеля-Дарбу) появились в основополагающих работах П.~Л.~Чебышёва. Отметим, что в 19-м веке роль цепных дробей в анализе была сравнима с ролью ряда Тейлора.

  Для диаграммы Дынкина (конечного графа без петель и кратных рёбер) $G_{\bar n}$ с множеством вершин $\bar n=\{1,\dots , n\}$ и весами на рёбрах $a_{ij}=a_{ji}$ (условие $a_{ij}=0$ равносильно отсутствию ребра $ij$ в графе) многочлен Кокстера (характеристический многочлен монодромии в случае особенностей \cite{2}) можно определить формулой 
$$
     G_{\bar n}^*(q)=\det {(qS+S^t)},
$$
где $S$ -- верхнетреугольная матрица с единицами на диагонали и элементами  $-a_{ij}$ выше диагонали \cite{5}. Мы будем работать с симметризованным многочленом
$$
   G_{\bar n}(q)=q^{-n}G_{\bar n}^*(q^2)=\det {(qS+q^{-1}S^t)}.
$$ 
В теории особенностей и теории узлов матрица $S$ отвечает форме Зейферта. Многочлен Александера-Конвея узла равен ${(-1)}^n\det {(qS-q^{-1}S^t)}$. Положим для всех $i$ $a_{ii}=-2$ и обозначим через $C=S+S^t$ матрицу $(-a_{ij})$. Для аффинных диаграмм Дынкина эта матрица совпадает с матрицей Картана. Для диаграмм Дынкина особенностей функций нечётного числа переменных матрица $C$ задаёт индекс пересечения в отмеченном базисе гомологий слоя Милнора \cite{2}. Матрица $qS+S^t$ задаёт  $q$-индекс пересечения \cite{1} и полученные с её помощью разложения в цепные дроби являются $q$-аналогами цепных дробей, появляющихся в теории разрешений особенностей и исчислении Кирби (возникает вопрос о существовании $q$-аналогов этих теорий -- возможно, такие обобщения могли бы основываться на результатах А.~Б.~Гивенталя о $q$-монодромии \cite{1}). Тождества Кристоффеля-Дарбу для рядов Пуанкаре групп Клейна и многочленов Кокстера могут оказаться полезными в обратном направлении -- интересно было бы выяснить смысл этих тождеств при $q=1$ в теории особенностей, теории узлов и теории алгебр Ли, отвечающих группам Клейна. В теории групп Клейна они представляют собой квадратичные соотношения между кратностями вхождения неприводимых представлений группы Клейна в ограничения на эту группу неприводимых представлений группы $SU(2)$. По-видимому, эти соотношения можно описать в терминах колец представлений группы $SU(2)$ и групп Клейна.

  В.~Эбелинг интерпретировал вектор рядов Пуанкаре каждой группы Клейна как решение линейной системы с матрицей $(z-2)E+C$, где $C$ -- матрица Картана соответствующей аффинной диаграммы Дынкина \cite{7}. Уравнения в этой линейной системе можно рассматривать как обобщения последовательных делений с остатком в алгоритме Евклида и, наряду с формулой дополнения по Шуру, использовать их для разложения в ветвящиеся цепные дроби и для доказательства аналогов тождества Кристоффеля-Дарбу. В \cite{7} использовался следующий известный факт: для графа, не имеющего циклов нечётной длины, характеристический многочлен $\det {((z-2)E+C)}$ совпадает с многочленом Кокстера. Используя это равенство и формулу для характеристического многочлена графа \cite{ 9}, мы покажем в п.~3, что для каждой из аффинных диаграмм Дынкина $\tilde A_n$, $\tilde D_n$, $\tilde E_6$, $\tilde E_7$, $\tilde E_8$ вектор рядов Пуанкаре соответствующей группы Клейна пропорционален вектору, элементы которого определяются многочленами Кокстера поддиаграмм. Коэффициент пропорциональности совпадает с умноженным на $q$ многочленом Кокстера соответствующей аффинной диаграммы Дынкина за исключением цикла нечётной длины $\tilde A_{2m}$. В этом случае он равен умноженному на $q$ характеристическому многочлену графа $\tilde A_{2m}$. В случае ряда Пуанкаре для инвариантов группы Клейна эти утверждения содержатся в \cite{7}.

  Для единой формулировки некоторых фактов о рядах Пуанкаре групп Клейна удобно считать, что для каждой аффинной диаграммы Дынкина существует ещё одна (не принадлежащая диаграмме) вершина, соединённая с аффинной вершиной ребром веса $1$, и ряд Пуанкаре для этой дополнительной вершины равен $q^{-1}$. Возникает вопрос, существуют ли интерпретации такой вершины в теории представлений групп Клейна или в теории клейновых особенностей?

  Результаты из \cite{9} и \cite{18} позволяют получить несколько фактов, связанных с корнями рядов Пуанкаре: разложение в сумму простейших дробей, представление безутиана и вронскиана из тождества Кристоффеля-Дарбу в виде симметрических функций от корней, выражение одного ряда Пуанкаре через корни другого, перемежаемость корней и связь с последовательностями Штурма, дискретная ортогональность, связанная с подстановкой корней в тождество Кристоффеля-Дарбу. Мы не будем подробно останавливаться на этих темах.

  Дж.~Левин описал перестройку зацепления в узел, при которой контролируемо меняется многочлен Александера-Конвея \cite{19}. А именно, отношение многочленов Александера-Конвея зацепления и узла определяется инвариантами Милнора  промежуточного объекта перестройки -- некоторого стринг-зацепления (обобщённой косы, в которой нитям разрешается не быть монотонными). В той же ситуации отношение многочленов Александера-Конвея можно описать в терминах $q$-индекса зацепления \cite{30}. Представляет интерес следующий вопрос: какие дополнительные структуры может сохранять перестройка Дж.~Левина? Мы частично ответим на этот вопрос в случаях, когда зацепление строится по плоской кривой согласно конструкции А'Кампо \cite{11} и когда зацепление строится по хордовой диаграмме согласно конструкции Хиронаки \cite{12}. Нас будет интересовать возможность разрезать отвечающее плоской кривой или хордовой диаграмме зацепление в нескольких точках и затем склеить образовавшиеся концы таким образом, что полученное зацепление также будет соответствовать некоторой плоской кривой или хордовой диаграмме. Другими словами, мы найдём перестройки зацеплений, реализуемые перестройками плоских кривых или хордовых диаграмм, в частности, мы реализуем перестройкой Левина разрезание плоской кривой в неособой точке и удаление хорды из хордовой диаграммы. Объект, полученный после разрезания, мы интерпретируем как стринг-зацепление. Это позволит в п.~4 применить теорему Левина и, в частности, выразить через инварианты Милнора стринг-зацеплений следующие отношения многочленов Кокстера евклидовых и аффинных диаграмм Дынкина  
$$
\frac{A_{2n+1}(q)}{A_{2n}(q)},\frac{A_{2n+1}(q)}{A_{2n+2}(q)},\frac{\tilde E_6(q)}{E_6(q)},\frac{D_5(q)}{E_6(q)},\frac{D_7(q)}{E_8(q)},\frac{E_7(q)}{E_8(q)}
$$
Пересечение этого списка со списком Примера 5 даёт отношения рядов Пуанкаре (или обратные к ним), которые выражаются через инварианты Милнора стринг-зацеплений.

  Для любого зацепления $K$ обозначим через $A_K(q)$ многочлен Александера-Конвея этого зацепления. Для стринг-зацепления $L$, как и для косы, определено замыкание $\tilde L$ и образ в представлении Бурау $\beta (L)$. В \cite{14} доказан следующий аналог формулы Левина, в которой роль перестройки играет умножение на любую обычную косу $B$,
$$
\frac{A_{\tilde L}(q)}{A_{\tilde {BL}}(q)}=\frac{\det {(E-\beta (L))}}{\det {(E-\beta (BL))}}.\eqno (1)
$$ 
Мы покажем в п. 4, что разрезание плоской кривой в неособой точке и удаление хорды из хордовой диаграммы реализуются умножением на стандартную образующую группы кос. В частности, можно выразить через представление Бурау отношения многочленов Кокстера евклидовых и аффинных диаграмм Дынкина из списка в предыдущем абзаце.

  Для неприводимой особенности кривой ряд Пуанкаре кольца функций на этой кривой совпадает с дзета-функцией монодромии (она определяется характеристическим многочленом монодромии) \cite{3}. Возможно, результаты о факторизации многочлена Александера-Конвея (на алгебраическом уровне эта факторизация реализуется формулой дополнения по Шуру) \cite{14}, \cite{19} и приведённые ниже формулы могут послужить основой для рекурсивного построения теории, изучающей связи дзета-функций монодромии и рядов Пуанкаре колец функций на особенностях. Сформулируем более конкретную гипотезу: отношение рядов Пуанкаре колец функций для близких (по примыканию или по расположению в серии) особенностей кривых определяется инвариантами Милнора стринг-зацепления, которое является промежуточным объектом при перестройке узла одной особенности в узел другой. Многочлен Кокстера евклидовой диаграммы Дынкина $A_n$ совпадает с $h$-многочленом $n$-симплекса или многочленом Пуанкаре соответствующего торического многообразия -- $n$-мерного проективного пространства. Можно надеятся, что существуют похожие связи между стринг-зацеплениями и отношениями $h$-многочленов простого многогранника и его граней или многочленов Пуанкаре соответствующих им торических многообразий.  В работах о связях дзета-функций и рядов Пуанкаре обычно используется представление этих объектов в виде произведений циклотомических многочленов, тогда как в работах \cite{14}, \cite{19} и в нашей статье существенно детерминантное представление. Это замечание подводит нас к следующей задаче, которую можно считать переформулировкой сказанного выше: перевести теорию о связях дзета-функций и рядов Пуанкаре с языка разрешений особенностей на язык форм Зейферта, в частности, найти прямое доказательство рекуррентных соотношений для рядов Пуанкаре колец функций на особенностях (они вытекают из рекуррентных соотношений для характеристического многочлена монодромии) и описать двойственность Саито для многочленов Кокстера \cite{7} в терминах связанных с формой Зейферта определителей и рекуррентных соотношений.

  \section {Дополнение по Шуру и многочлены Кокстера}

  Обозначим $q+q^{-1}$ через $z$. Используя в качестве $A_{11}$ матрицу $(z)$, расположенную на пересечении первых строки и столбца матрицы $qS+q^{-1}S^t$, мы можем переписать формулу дополнения по Шуру следующим образом
$$
 G_{\bar n}(q)=zG_{\bar n \setminus 1}(q)-\sum_{2\le i\le n}{a_{1i}^2G_{\bar n \setminus {\{1,i\}}}(q)}-\sum_{2\le i \ne j\le n}{a_{1i}a_{1j}P_{ij}(q)},\eqno (2)
$$
где $P_{ij}(q)$ $(=P_{ji}(q^{-1}))$ -- алгебраическое дополнение элемента $-a_{ij}$, $2\le i \ne j\le n$, в матрице, полученной из матрицы $qS+q^{-1}S^t$ вычёркиванием первых строки и столбца. В частности, эта формула позволяет одновременно с вычислением определителя перейти от переменной  $q$ к переменной $z=q+q^{-1}$. В теории узлов обычно используется матрица $-qS+q^{-1}S^t$ и переход от переменной  $q$ к переменной $t=q-q^{-1}$ является переходом к многочлену Конвея. Переход к переменной $z=q+q^{-1}$ может быть полезен, например, по следующей причине: для характеристического многочлена монодромии мы получим многочлен с вещественными корнями, к которому можно применить хорошо развитую теорию таких многочленов.

  Следующий частный (неоднократно переоткрывавшийся) случай формулы~(2) использовался при изучении многочленов Кокстера в контексте теории особенностей, теории алгебраических чисел, теории колчанов, теории представлений ассоциативных алгебр и частично упорядоченных множеств (ссылки можно найти в \cite{28}). Джойном диаграмм Дынкина $T_i$ с отмеченными вершинами $v_i$ называется диаграмма Дынкина $T_{\bar n}$, полученная добавлением к объединению диаграмм $T_i$ ещё одной вершины $v$ и рёбер веса $1$, соединяющих вершину $v$ со всеми отмеченными вершинами. Из формулы (1) вытекает следующее равенство (через ${\bar T}_i$ обозначим диаграмму Дынкина, полученную из диаграммы $T_i$ удалением вершины $v_i$). 
$$
T_{\bar n}(q)=\prod_i {T_i(q)}(z-\sum_j{\frac{{\bar T}_j (q)}{T_j(q)}}).\eqno (3)
$$

  Формулы~(2) и (3) можно использовать для разложений в цепные дроби
$$
 \frac{G_{\bar n \setminus 1}(q)}{G_{\bar n}(q)}
=\frac{1}{z-\sum_{2\le i\le n}{a_{1i}^2\frac{G_{\bar n \setminus {\{1,i\}}}(q)}{G_{\bar n \setminus 1}(q)}}-\sum_{2\le i \ne j\le n}{a_{1i}a_{1j}\frac{P_{ij}(q)}{G_{\bar n\setminus 1}(q)}}},
$$
$$
\frac{\prod_i {T_i(q)}}{T_{\bar n}(q)}=\frac{1}{z-\sum_j{\frac{{\bar T}_j (q)}{T_j(q)}}},
$$

 Пример 1. Выпишем цепные дроби для аффинных диаграмм Дынкина $\tilde D_n$, $\tilde E_6$, $\tilde E_7$, $\tilde E_8$. Количество появлений переменной $z$ в каждой цепной дроби равно числу вершин соответствующей аффинной диаграммы Дынкина. Форма цепной дроби в определённом смысле совпадает с формой соответствующей диаграммы. В знаменателях цепной дроби, содержащих больше двух слагаемых, первое слагаемое $z$ соответствует вершине ветвления в диаграмме Дынкина. Например, в диаграмме $\tilde D_n$ одна вершина ветвления при $n=4$ и две при $n>4$. Поэтому для диаграммы $\tilde D_n$ мы выделяем два случая. Из формулы Эбелинга \cite{7} следует, что для каждой диаграммы цепная дробь равна умноженному на $q$ ряду Пуанкаре для инвариантов соответствующей группы Клейна.
$$
  \frac{D_4(q)}{\tilde D_4(q)}=
      \cfrac{1}{
      z-\cfrac{1}{
       z-\cfrac{1}{z}-\cfrac{1}{z}-\cfrac{1}{z}}}
$$
$$
  \frac{D_n(q)}{\tilde D_n(q)}=
      \cfrac{1}{
      z-\cfrac{1}{
       z-\cfrac{1}{
        z-\cfrac{1}{
        z-\dots-\cfrac{1}{
         z-\cfrac{1}{
          z-\cfrac{1}{z}-\cfrac{1}{z}}}}}-\cfrac{1}{z}}}
$$
$$
  \frac{E_6(q)}{\tilde E_6(q)}=
      \cfrac{1}{
      z-\cfrac{1}{
       z-\cfrac{1}{
        z-\cfrac{1}{z-\cfrac{1}{z}}-\cfrac{1}{z-\cfrac{1}{z}}}}}
$$ 
$$
  \frac{E_7(q)}{\tilde E_7(q)}=
        \cfrac{1}{
      z-\cfrac{1}{
       z-\cfrac{1}{
        z-\cfrac{1}{
         z-\cfrac{1}{
          z-\cfrac{1}{
           z-\cfrac{1}{z}}}-\cfrac{1}{z}}}}}
$$
$$
  \frac{E_8(q)}{\tilde E_8(q)}=
        \cfrac{1}{
      z-\cfrac{1}{
       z-\cfrac{1}{
        z-\cfrac{1}{
         z-\cfrac{1}{
          z-\cfrac{1}{
           z-\cfrac{1}{
            z-\cfrac{1}{z}}-\cfrac{1}{z}}}}}}}
$$

  Безутианом и вронскианом функций $f$ и $g$ называются выражения
$$
  Bez(f,g)=\frac{f(x)g(y)-f(y)g(x)}{x-y},
$$
$$
  Wr(f,g)=\lim_{y\to x}{Bez(f,g)}=f'g-fg'.
$$
Для многочленов $f$ и $g$ вронскиан также может быть представлен следующей формулой
$$
  Wr(f,g)=m^{-1}J(F,G)|_{y=1},
$$
где $F=F(x,y)$ и $G=G(x,y)$ -- однородные формы степени $m$, $F(x,1)=f$, $G(x,1)=g$ и $J(F,G)$ -- якобиан форм $F$ и $G$ \cite{18}.  
\begin{prop}\label{prop1} Для произвольной диаграммы Дынкина $G_{\bar n}$ справедливы соотношения
$$
 Bez(G_{\bar n}(q),G_{\bar n \setminus 1}(q))=(1-(xy)^{-1})G_{\bar n \setminus 1}(x)G_{\bar n \setminus 1}(y)\eqno (4)
$$
$$
+\sum_{2\le i\le n}{a_{1i}^2Bez(G_{\bar n \setminus 1}(q), G_{\bar n \setminus {\{1,i\}}}(q))}+\sum_{2\le i \ne j\le n}{a_{1i}a_{1j}Bez(G_{\bar n \setminus 1}(q), P_{ij}(q))},
$$
$$
 Wr(G_{\bar n}(q),G_{\bar n \setminus 1}(q))=(1-x^{-2})G_{\bar n \setminus 1}^2(x)\eqno (5)
$$
$$
+\sum_{2\le i\le n}{a_{1i}^2Wr(G_{\bar n \setminus 1}(q), G_{\bar n \setminus {\{1,i\}}}(q))}+\sum_{2\le i \ne j\le n}{a_{1i}a_{1j}Wr(G_{\bar n \setminus 1}, P_{ij}(q))},
$$
в частности,
$$
 Bez(T_{\bar n}(q),\prod_i {T_i(q)})=(1-(xy)^{-1})\prod_i {T_i(x)T_i(y)}+\sum_j{Bez(T_j(q),{\bar T}_j(q))},
$$
$$
Wr(T_{\bar n}(q),\prod_i {T_i(q)})=(1-x^{-2})\prod_i {T_i^2(x)}+\sum_j{Wr(T_j(q),{\bar T}_j(q))}.
$$
\end{prop}
 
Доказательство. Формула для вронскиана получается предельным переходом из формулы для безутиана. Докажем формулу (4) по аналогии с доказательством тождества Кристоффеля-Дарбу в \cite{18}.
$$
\frac{G_{\bar n}(x)G_{\bar n \setminus 1}(y)-G_{\bar n}(y)G_{\bar n \setminus 1}(x)}{G_{\bar n \setminus 1}(x)G_{\bar n \setminus 1}(y)}
=\frac {G_{\bar n}(x)}{G_{\bar n \setminus 1}(x)}
-\frac {G_{\bar n}(y)}{G_{\bar n \setminus 1}(y)}
$$
$$
=x+x^{-1}-\sum_{2\le i\le n}{a_{1i}^2\frac{G_{\bar n \setminus {\{1,i\}}}(x)}{G_{\bar n \setminus 1}(x)}}-\sum_{2\le i \ne j\le n}{a_{1i}a_{1j}\frac{P_{ij}(x)}{G_{\bar n\setminus 1}(x)}}
$$
$$
-y-y^{-1}+\sum_{2\le i\le n}{a_{1i}^2\frac{G_{\bar n \setminus {\{1,i\}}}(y)}{G_{\bar n \setminus 1}(y)}}+\sum_{2\le i \ne j\le n}{a_{1i}a_{1j}\frac{P_{ij}(y)}{G_{\bar n\setminus 1}(y)}}
$$
$$
=x-y+x^{-1}-y^{-1}+\sum_{2\le i\le n}{a_{1i}^2\frac{G_{\bar n\setminus 1}(x)G_{\bar n \setminus {\{1,i\}}}(y)-G_{\bar n\setminus 1}(y)G_{\bar n \setminus {\{1,i\}}}(x)}{G_{\bar n \setminus 1}(x)G_{\bar n \setminus 1}(y)}}
$$
$$
+\sum_{2\le i \ne j\le n}{a_{1i}a_{1j}\frac{G_{\bar n\setminus 1}(x)P_{ij}(y)-G_{\bar n\setminus 1}(y)P_{ij}(x)}{G_{\bar n \setminus 1}(x)G_{\bar n \setminus 1}(y)}}
$$
и остаётся домножить на $(G_{\bar n \setminus 1}(x)G_{\bar n \setminus 1}(y))/(x-y)$.

  Пример 2. Пусть в диаграмме Дынкина $G_{\bar n}$ поддиаграмма на вершинах $1,\dots ,k$ совпадает с евклидовой диаграммой Дынкина $A_k$ и вершины $1,\dots ,k-1$ не связаны рёбрами с поддиаграммой $G_{{\bar n}\setminus {\bar k}}$. Нумерация вершин выбрана так, что $a_{12}=a_{23}=\dots a_{k-1k}=1$. Обозначим через $C_i$ поддиаграмму $G_{{\bar n}\setminus {\bar i}}$, $1\le i\le k$, $C_0=G_{\bar n}$. Тогда имеем соотношения
$$
  C_{i-1}(q)-zC_i(q)+C_{i+1}(q)=0, 1\le i\le {k-1},
$$
которые можно представить в следующем виде
$$
\begin{pmatrix}
C_{i-1}(q)&C_i(q)\\
C_i(q)&C_{i+1}(q)
\end{pmatrix}=
\begin{pmatrix}
0&1\\
-1&z
\end{pmatrix}^{i-1}
\begin{pmatrix}
C_0(q)&C_1(q)\\
C_1(q)&C_2(q)
\end{pmatrix},
$$
Используя эти соотношения и суммируя равенства вида
$$
\frac{C_{i-1}(q)}{C_i(q)}-\frac{C_i(q)}{C_{i+1}(q)}=
\begin{vmatrix}
C_{i-1}(q)&C_i(q)\\
C_i(q)&C_{i+1}(q)
\end{vmatrix}
\frac{1}{C_i(q)C_{i+1}(q)},
$$
получим для $1\le i\le {k-1}$ следующие тождества
$$
\frac{C_{i-1}(q)}{C_i(q)}=
\begin{vmatrix}
C_0(q)&C_1(q)\\
C_1(q)&C_2(q)
\end{vmatrix}
\sum_{j=i+1}^k{\frac{1}{C_{j-1}(q)C_{j}(q)}}+\frac{C_{k-1}(q)}{C_k(q)}.
$$
Разложения в цепные дроби и формулы для безутиана и вронскиана принимают следующий вид (количество появлений переменной $z$ в цепной дроби равно $k-i$)
$$
    \frac{C_i(q)}{C_{i-1}(q)}=
     \cfrac{1}{
      z-\cfrac{1}{
       z-\dots-\cfrac{1}{
        z-\cfrac{1}{
         z-{\frac{C_k(q)}{C_{k-1}(q)}}}}}},
$$
$$
  Bez(C_{i-1}(q),C_i(q))=(1-(xy)^{-1})\sum_{j=i}^{k-1}{C_j(x)C_j(y)}
  +Bez(C_{k-1}(q),C_k(q)),
$$
$$
   Wr(C_{i-1}(q),C_i(q))=(1-x^{-2})
\sum_{j=i}^{k-1}{C_j^2(x)}
+Wr(C_{k-1}(q),C_k(q)).
$$

  Характеристическим многочленом графа $G_{\bar n}$ называется следующий определитель
$$
     G_{\bar n}^\#(z)=\det {((z-2)E+C)}.
$$
Приведённые выше формулы для многочлена Кокстера остаются справедливыми для характеристического многочлена, если заменить миноры матрицы $qS+q^{-1}S^t$ на соответствующие миноры матрицы $(z-2)E+C$. Обозначим через $H_{ij}(z)(=H_{ji}(z))$ алгебраические дополнения элементов матрицы $(z-2)E+C$, через $Q_{ij}$ -- множество путей без повторения вершин, идущих в графе $G_{\bar n}$ из вершины $i$ в вершину $j$, через $d_{ij}^k$ -- сумму весов путей в графе $G_{\bar n}$, имеющих длину $k$ и идущих из вершины $i$ в вершину $j$ (весом $a(P)$ пути $P$ называется произведение весов входящих в него рёбер) Следующие тождества (а также ссылки на первоисточники) имеются в \cite{9}:
$$
H_{ij}(z)=\sum_{Q\in {Q_{ij}}}{a(Q)(G_{\bar n}^\#\setminus Q)(z)},\eqno (6)
$$
$$
H_{ij}^2(z)=
\begin{vmatrix}
G_{\bar n\setminus i}^\#(z)&G_{\bar n\setminus {\{i,j\}}}^\#(z)\\
G_{\bar n}^\#(z)&G_{\bar n\setminus j}^\#(z)
\end{vmatrix},i\ne j,\eqno (7)
$$
$$
Bez(G_{\bar n}^\#(z),H_{ij}(z))=\sum_{k\in {\bar n}}{H_{ik}(x)H_{jk}(y)},\eqno (8)
$$
$$
Wr(G_{\bar n}^\#(z),H_{ij}(z))=\sum_{k\in {\bar n}}{H_{ik}^2(x)},\eqno (9)
$$
$$
\frac{H_{ij}(z)}{G_{\bar n}^\#(z)}=\sum_{0\le k}{d_{ij}^kz^{-k-1}}.\eqno (10)
$$

 Если граф является деревом, то при любой нумерации вершин его характеристический многочлен совпадает с многочленом Кокстера. Этот известный факт вытекает, например, из формуы (3) по индукции. Известен также более широкий класс графов со специальной нумерацией вершин, для которых $G_{\bar n}^\#(z)=G_{\bar n}(q)$. Эти графы называются двудольными и характеризуются отсутствием циклов нечётной длины. По определению, вершины такого графа можно разбить на два множества, в каждом из которых нет соединённых рёбрами вершин. Нумеруются сначала вершины одного множества, а затем другого. Факт совпадения характеристического многочлена и многочлена Кокстера для двудольного графа с такой нумерацией вершин является прямым следствием формулы дополнения по Шуру. Это совпадение использовалось при доказательстве формулы Эбелинга \cite{7}. Мы используем его для обобщения формулы Эбелинга. Среди диаграмм Дынкина особенностей деревьями являются только стандартные диаграммы Дынкина простых особенностей $A_n$, $D_n$, $E_6$, $E_7$, $E_8$. В \cite{6} имеются примеры диаграмм Дынкина не простых особенностей, которые не имеют циклов нечётной длины. Отметим ещё одно следствие совпадения характеристического многочлена и многочлена Кокстера. Подграф двудольного графа является двудольным и диагональные миноры матрицы $(z-2)E+C$ равны многочленам Кокстера соответствующих подграфов. Поэтому тождества для диагональных миноров симметрической матрицы \cite{13} дают соотношения для многочленов Кокстера подграфов.

 Пример 3. Рассмотрим аффинные диаграмм Дынкина $\tilde A_n$, $\tilde D_n$, $\tilde E_6$, $\tilde E_7$, $\tilde E_8$ с нумерацией вершин, в которой вершина с номером $0$ совпадает с вершиной, добавляемой при переходе от евклидовой диаграммы Дынкина к аффинной. Веса всех рёбер равны $1$. Используя факты из предыдущего абзаца и формулу (6), для каждой диаграммы выразим многочлены $H_{i0}(z)$ через многочлены Кокстера поддиаграмм. Для аффинной диаграммы Дынкина $\tilde A_n$ имеем (полагаем $A_0(z)=1$) 
$$
H_{00}(z)=A_n(q),H_{j0}(z)=A_{i-2}(q)+A_{n-i+1}(q),1\le j\le {n}.
$$
Напомним, что аффинная диаграмма Дынкина $\tilde A_n$ является циклом с $n+1$ вершиной и поэтому из одной вершины в другую можно попасть двумя путями. Оставшиеся многочлены $H_{j0}(z)$ представлены на диаграммах. Многочлен $H_{00}(z)$ всегда совпадает с многочленом Кокстера соответствующей евклидовой диаграммы Дынкина. В п.~3 мы увидим, что для каждой из диаграмм $\tilde A_n$, $\tilde D_n$, $\tilde E_6$, $\tilde E_7$, $\tilde E_8$ набор многочленов $H_{i0}(z)$ пропорционален набору рядов Пуанкаре соответствующей группы Клейна. Коэффициент пропорциональности совпадает с умноженным на $q$ многочленом Кокстера соответствующей аффинной диаграммы Дынкина за исключением цикла нечётной длины $\tilde A_{2m}$. В этом случае он равен умноженному на $q$ характеристическому многочлену графа $\tilde A_{2m}$.

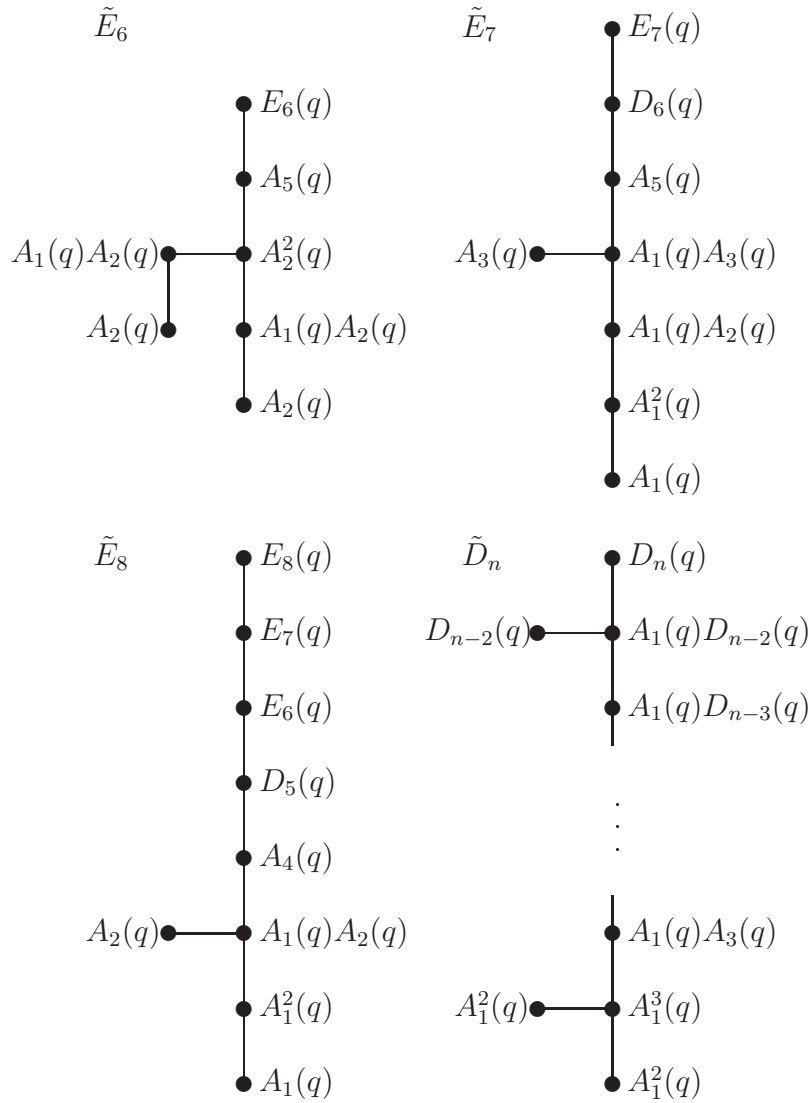
\begin{picture}(120,70)
\put(11,64){$\tilde E_6$}
\put(31,15){\line(0,1){40}}
\put(31,35){\line(-1,0){10}}
\put(21,35){\line(0,-1){10}}
\put(21,35){\circle*{2}}
\put(21,25){\circle*{2}}
\multiput(31,15)(0,10){5}{\circle*{2}}
\put(33,14){$A_2(q)$}
\put(33,24){$A_1(q)A_2(q)$}
\put(33,34){$A_2^2(q)$}
\put(33,44){$A_5(q)$}
\put(33,54){$E_6(q)$}
\put(10,24){$A_2(q)$}
\put(0,34){$A_1(q)A_2(q)$}
\put(60,64){$\tilde E_7$}
\put(80,5){\line(0,1){60}}
\put(80,35){\line(-1,0){10}}
\put(70,35){\circle*{2}}
\multiput(80,5)(0,10){7}{\circle*{2}}
\put(82,4){$A_1(q)$}
\put(82,14){$A_1^2(q)$}
\put(82,24){$A_1(q)A_2(q)$}
\put(82,34){$A_1(q)A_3(q)$}
\put(82,44){$A_5(q)$}
\put(82,54){$D_6(q)$}
\put(82,64){$E_7(q)$}
\put(59,34){$A_3(q)$}
\end{picture}
 
\begin{picture}(120,80)
\put(11,74){$\tilde E_8$}
\put(31,5){\line(0,1){70}}
\put(31,25){\line(-1,0){10}}
\put(21,25){\circle*{2}}
\multiput(31,5)(0,10){8}{\circle*{2}}
\put(33,4){$A_1(q)$}
\put(33,14){$A_1^2(q)$}
\put(33,24){$A_1(q)A_2(q)$}
\put(33,34){$A_4(q)$}
\put(33,44){$D_5(q)$}
\put(33,54){$E_6(q)$}
\put(33,64){$E_7(q)$}
\put(33,74){$E_8(q)$}
\put(10,24){$A_2(q)$}
\put(60,74){$\tilde D_n$}
\put(80,5){\line(0,1){25}}
\put(80,50){\line(0,1){25}}
\put(80,15){\line(-1,0){10}}
\put(70,15){\circle*{2}}
\multiput(80,5)(0,10){3}{\circle*{2}}
\multiput(80,55)(0,10){3}{\circle*{2}}
\multiput(80,36)(0,3){3}{.}
\put(82,4){$A_1^2(q)$}
\put(82,14){$A_1^3(q)$}
\put(82,24){$A_1(q)A_3(q)$}
\put(82,54){$A_1(q)D_{n-3}(q)$}
\put(82,64){$A_1(q)D_{n-2}(q)$}
\put(80,65){\line(-1,0){10}}
\put(70,65){\circle*{2}}
\put(82,74){$D_n(q)$}
\put(59,14){$A_1^2(q)$}
\put(55,64){$D_{n-2}(q)$}
\end{picture}

 Пример 4. Следующее равенство является аналогом тождества (1) для аффинной диаграммы Дынкина $\tilde A_n$
$$
\frac{H_{11}^\#(z)}{\tilde A_n^\#(z)}=\frac{A_n^\#(z)}{\tilde A_n^\#(z)}=\frac{1}{z-2\frac{A_{n-1}^\#(z)+1}{A_n^\#(z)}}.
$$
Используя соотношение
$$
A_{k+1}^\#(z)=zA_k^\#(z)-A_{k-1}^\#(z),
$$
мы можем продолжить разложение в цепную дробь. Из формулы Эбелинга \cite{7} следует, что результат разложения совпадает с умноженным на $q$ рядом Пуанкаре для инвариантов соответствующей диаграмме $\tilde A_n$ группы Клейна.
$$
\frac{A_n^\#(z)}{\tilde A_n^\#(z)}=
\cfrac{1}{
      z-\cfrac{1}{
       z-\dots-\cfrac{1}{
        z-\cfrac{1}{r}}}-\cfrac{1}{
       z-\dots-\cfrac{1}{
        z-\cfrac{1}{r}}}},
$$
где для нечётного $n$ $r=z/2$, для чётного $n$ $r=1$, для нечётного $n$ переменная $z$ появляется в цепной дроби $n$ раз, а для чётного $n$ -- $n+1$ раз.

 Известен способ обобщения тождества Кристоффеля-Дарбу \cite{23}. Из формулы (2) для наборов переменных $x_1,\dots ,x_m$ и $y_1,\dots ,y_k$ вытекает матричное равенство ($1\le l\le m$, $1\le s\le k$, в скобках после безутиана мы выписываем переменные, от которых он зависит)
$$
(Bez(G_{\bar n}^\#(z),H_{ij}(z))(x_l,y_s))=(H_{ik}(x_l))(H_{jk}(y_s))^t.\eqno (11)
$$
Если $m=k\le n$, то мы можем применить к этому равенству тождество Бине-Коши:  
$$
\det {(Bez(G_{\bar n}^\#(z),H_{ij}(z))(x_l,y_s))}=\sum \det {(H_{ik_r}(x_l))}\det {(H_{jk_r}(y_s))},\eqno (12)
$$
где суммирование ведётся по всем $m$-подмножествам $1\le k_1<\dots <k_m\le n$ в множестве $\bar n$.

  \section {Ряды Пуанкаре и многочлены Кокстера}

  Для каждого неприводимых представления $r_i$ конечной подгруппы $B\subset SU(2)$ (бинарной полиэдральной группы) Б.~Костант определил ряды Пуанкаре
$$
        P_i(q)=\sum_{0\le n} m_{n,i}q^n,
$$
где $m_{n,i}$ -- кратность представления $r_i$ в ограничении на $B$ $n$-й симметрической степени стандартного двумерного представления группы $SU(2)$ (известно, что симметрические степени стандартного представления исчерпывают неприводимые представления группы $SU(2)$) \cite{16}, \cite{17}. Для этих рядов он получил следующие формулы
$$
     P_i(q)=\frac{Z_i(q)}{(1-q^a)(1-q^b)},\eqno (13)
$$
где $Z_i(q)$ -- некоторые многочлены с целыми неотрицательными коэффициентами, степени которых не превосходят $h$ -- суммы размерностей неприводимых представлений группы $B$, $a$ и $b$ можно найти из условий $a+b=h+2$ и $ab=2\vert B\vert$. Отметим арифметическое следствие этого результатата Б.~Костанта: для любой группы Клейна $B$ натуральное число $(h+2)^2-8\vert B\vert$ является полным квадратом,так как натуральные числа $a$ и $b$ являются корнями квадратного уравнения $x^2-(h+2)x+2\vert B\vert=0$. Ниже мы получим некоторые формулы для отношений рядов Пуанкаре $P_i(q)$, очевидно эти отношения совпадают с отношениями соответствующих многочленов $Z_i(q)$. Характеры группы $B$ выражаются через специализации многочленов $Z_i(q)$ \cite{24} и поэтому формулы для отношений многочленов $Z_i(q)$ дают соотношения для характеров.

  Соответствие Маккея сопоставляет группе $B$ аффинную диаграмму Дынкина некоторой алгебры Ли. Вершины диаграммы взаимно однозначно соответствуют неприводимым представлениям группы $B$, причём, тождественное представление $r_0$ соответствует вершине, добавляемой при переходе от евклидовой диаграммы Дынкина к аффинной (назовём эту вершину аффинной). Б.~Костант определил коэффициенты многочленов $Z_i(q)$ в терминах орбит действия элемента Кокстера диаграммы на соответствующей системе корней, число $h$ равно порядку элемента Кокстера. Другие подходы к рядам Пуанкаре $P_i(q)$ и многочленам $Z_i(q)$ можно найти в \cite{10}, \cite{15}, \cite{20}, \cite{21}, \cite{27}. 

  Для аффинной диаграммы Дынкина $\tilde A_n$ $a=2$, $b=n+1$ и $Z_i(q)=q^i+q^{n-i+1}$, $i=0, 1, \dots, n$. Для аффинных диаграмм Дынкина $\tilde D_n$, $\tilde E_6$, $\tilde E_7$, $\tilde E_8$ числа $a$, $b$ и многочлены $Z_i(q)$ представлены на диаграммах (вектору $(m, n, \dots)$ отвечает многочлен $q^m+q^n+ \dots$)

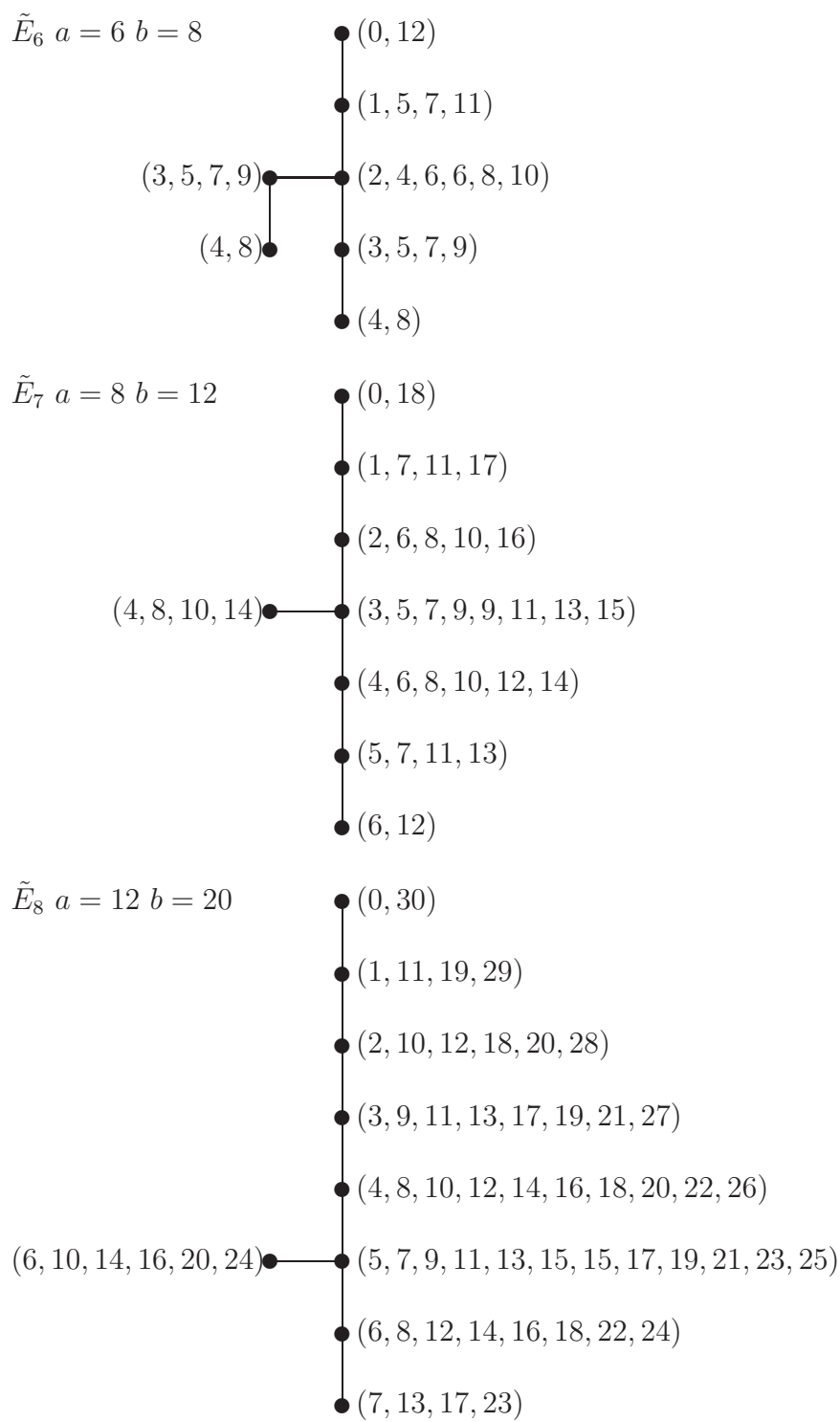
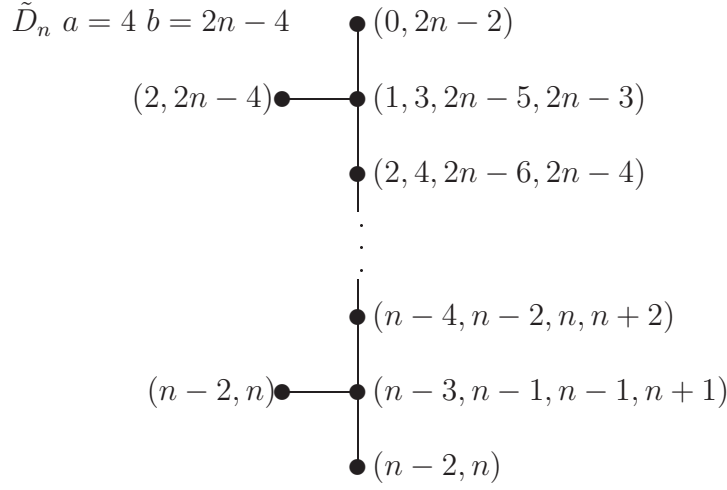
\begin{picture}(120,50)
\put(9,44){$\tilde E_6$ $a=6$ $b=8$}
\put(55,5){\line(0,1){40}}
\put(55,25){\line(-1,0){10}}
\put(45,25){\line(0,-1){10}}
\put(45,25){\circle*{2}}
\put(45,15){\circle*{2}}
\multiput(55,5)(0,10){5}{\circle*{2}}
\put(57,4){$(4, 8)$}
\put(57,14){$(3, 5, 7, 9)$}
\put(57,24){$(2, 4, 6, 6, 8, 10)$}
\put(57,34){$(1, 5, 7, 11)$}
\put(57,44){$(0, 12)$}
\put(35,14){$(4, 8)$}
\put(27,24){$(3, 5, 7, 9)$}
\end{picture}

\begin{picture}(120,70)
\put(9,64){$\tilde E_7$ $a=8$ $b=12$}
\put(55,5){\line(0,1){60}}
\put(55,35){\line(-1,0){10}}
\put(45,35){\circle*{2}}
\multiput(55,5)(0,10){7}{\circle*{2}}
\put(57,4){$(6, 12)$}
\put(57,14){$(5, 7, 11, 13)$}
\put(57,24){$(4, 6, 8, 10, 12, 14)$}
\put(57,34){$(3, 5, 7, 9, 9, 11, 13, 15)$}
\put(57,44){$(2, 6, 8, 10, 16)$}
\put(57,54){$(1, 7, 11, 17)$}
\put(57,64){$(0, 18)$}
\put(23,34){$(4, 8, 10, 14)$}
\end{picture}
 
\begin{picture}(120,80)
\put(9,74){$\tilde E_8$ $a=12$ $b=20$}
\put(55,5){\line(0,1){70}}
\put(55,25){\line(-1,0){10}}
\put(45,25){\circle*{2}}
\multiput(55,5)(0,10){8}{\circle*{2}}
\put(57,4){$(7, 13, 17, 23)$}
\put(57,14){$(6, 8, 12, 14, 16, 18, 22, 24)$}
\put(57,24){$(5, 7, 9, 11, 13, 15, 15, 17, 19, 21, 23, 25)$}
\put(57,34){$(4, 8, 10, 12, 14, 16, 18, 20, 22, 26)$}
\put(57,44){$(3, 9, 11, 13, 17, 19, 21, 27)$}
\put(57,54){$(2, 10, 12, 18, 20, 28)$}
\put(57,64){$(1, 11, 19, 29)$}
\put(57,74){$(0,30)$}
\put(9,24){$(6, 10, 14, 16, 20, 24)$}
\end{picture}

\begin{picture}(120,70)
\put(9,63){$\tilde D_n$ $a=4$ $b=2n-4$}
\put(55,5){\line(0,1){25}}
\put(55,39){\line(0,1){25}}
\put(55,15){\line(-1,0){10}}
\put(45,15){\circle*{2}}
\multiput(55,5)(0,10){3}{\circle*{2}}
\multiput(55,44)(0,10){3}{\circle*{2}}
\multiput(55,31)(0,3){3}{.}
\put(57,4){$(n-2, n)$}
\put(57,14){$(n-3, n-1, n-1, n+1)$}
\put(57,24){$(n-4, n-2, n, n+2)$}
\put(57,43){$(2, 4, 2n-6, 2n-4)$}
\put(57,53){$(1, 3, 2n-5, 2n-3)$}
\put(55,54){\line(-1,0){10}}
\put(45,54){\circle*{2}}
\put(57,63){$(0, 2n-2)$}
\put(27,14){$(n-2, n)$}
\put(25,53){$(2, 2n-4)$}
\end{picture}

   Нам удобно будет считать, что для каждой аффинной диаграммы Дынкина $T$ существует дополнительная (не принадлежащая $T$) вершина с номером $-1$, соединённая с имеющей номер $0$ аффинной вершиной ребром веса $1$, и $P_{-1}(q)=q^{-1}$, $Z_{-1}(q)=q^{-1}(1-q^a)(1-q^b)$, $H_{-10}(z)=T^\#(z)$ (для диаграмм $\tilde A_{2m-1}$, $\tilde D_n$, $\tilde E_6$, $\tilde E_7$, $\tilde E_8$ $T^\#(z)=T(q)$). Для  каждой вершины $i$ в диаграмме $T$ обозначим через $i{\to }$ вершину, которая соединена ребром с вершиной $i$ и расположена в диаграмме ближе к вершине $-1$, чем вершина $i$, в частности, вершина $0{\to }$ совпадает с вершиной $-1$. Если (в случае диаграммы $\tilde A_{2m+1}$) таких вершин окажется две, то берём любую из них. Для  каждой вершины $i$ обозначим через $Q_i$ множество вершин, кратчайшие пути из которых в вершину $-1$ содержат вершину $i$. Нумерация вершин в аффинной диаграмме $\tilde A_n$ соответствует обходу цикла.
 
  Используя формулу Клебша-Гордана и соответствие Маккея, В.~Эбелинг показал, что 
$$
  (P_0(q),P_1(q), \dots)((z-2)E+C)=(P_{-1}(q),0, \dots),\eqno (14)
$$
где C -- матрица Картана аффинной диаграммы Дынкина \cite{7} (некоторые дополнительные аспекты этой формулы представлены в \cite{29}). Из формул~(13) и (14) вытекает, что
$$
  (Z_0(q),Z_1(q), \dots)(zE+(C-2E))=(Z_{-1}(q),0, \dots).\eqno (15)
$$
Согласно определению и правилу Крамера, для любого графа имеем 
$$
  (H_{00}(z),H_{10}(z), \dots)((z-2)E+C)=(H_{-10}(z),0, \dots).\eqno (16)
$$
  Следующие утверждения вытекают из формул~(6)-(10) и (14)-(16). Равенства~(20) и (21) проще вывести из формул~(4) и (5). Упрощению формул способствует кососимметричность безутиана и вронскиана
$$
      Bez(f,g)=-Bez(g,f), Wr(f,g)=-Wr(g,f).
$$
Если многократным применением формул~(4) и (5) мы раскладываем безутиан или вронскиан в двух местах диаграммы, идя по циклу $\tilde A_{n}$ в разных направлениях и удаляясь от аффинной вершины, то в месте встречи этих разложений происходит сокращение промежуточных безутианов или вронскианов и остаются только слагаемые вида $(1-(xy)^{-1})P_k(x)P_k(y)$ или $(1-x^{-2})P_k^2(x)$.

\begin{prop}\label{prop2} 
  1) Для $-1\le i,j\le n$ имеем следующие равенства для рядов Пуанкаре группы Клейна, соответствующей любой из аффинных диаграмм Дынкина $\tilde A_n$, $\tilde D_n$, $\tilde E_6$, $\tilde E_7$, $\tilde E_8$,
$$   
\frac{P_i(q)}{P_{j}(q)}=\frac{Z_i(q)}{Z_j(q)}=\frac{H_{i0}(z)}{H_{j0}(z)},\eqno (17)
$$
в частности, для $i\ge 0$
$$
  P_i(q)=\frac{H_{i0}(z)}{qH_{-10}(z)}.
$$
Многочлены $H_{i0}(z)$ выражены через многочлены Кокстера в Примере~3. Исключение составляют многочлены $H_{-10}(z)$ для диаграмм $\tilde A_{2m}$ -- они совпадают с характеристическими многочленами графов $\tilde A_{2m}$ (полагаем $\tilde A_0(z)=z=A_1(z)$)
$$
q\tilde A_{2m}^\#(z) =q(z^{2m+1}-\sum_{i=1}^m{C_{2m+1}^i\tilde A_{2m-2i}^\#(z)}-2)=q^{-2m}(q^{2m+1}-1)^2.
$$
Для многочленов $H_{i0}(z)$, отвечающих каждой ветви в диаграммах $\tilde D_n$, $\tilde E_6$, $\tilde E_7$, $\tilde E_8$, справедливы формулы из примера~2.

  2) Для $1\le i\le n$ 
$$
P_i^2(q)=\frac{1}{qH_{-10}(z)}
\begin{vmatrix}
P_0(q)&T_i(q)\\
q^{-1}&\tilde T_i(q)
\end{vmatrix},
$$
где $T_i(q)$($\tilde T_i(q)$) -- многочлен Кокстера диаграммы, полученной удалением вершины $i$ из соответствующей евклидовой (аффинной) диаграммы Дынкина.

  3) Для $0\le i\le n$ имеем следующие равенства для аффинных диаграмм Дынкина $\tilde D_n$, $\tilde E_6$, $\tilde E_7$, $\tilde E_8$
$$ 
      Bez(P_{i_{\to }}(q),P_i(q))=(1-(xy)^{-1})\sum_{k\in Q_i}{P_k(x)P_k(y)},
$$
$$
      Wr(P_{i_{\to }}(q),P_i(q))=(1-x^{-2})\sum_{k\in Q_i}{P_k^2(x)},   
$$
в частности,
$$ 
      Bez(q^{-1},P_0(q))=(1-(xy)^{-1})\sum_{k\ge 0}{P_k(x)P_k(y)},\eqno (18)
$$
$$
      Wr(q^{-1},P_0(q))=(1-x^{-2})\sum_{k\ge 0}{P_k^2(x)},\eqno (19)   
$$
Для аффинной диаграммы Дынкина $\tilde A_n$ справедливы формулы~(18) и (19). Аналогами остальных формул являются следующие равенства, имеющие место при $1\le i\le j\le {n}$ (полагаем, что $P_{n+1}(q)=P_0(q)$),
$$       Bez(P_{i-1}(q),P_i(q))+Bez(P_j(q),P_{j+1}(q))=(1-(xy)^{-1})\sum_{k=i}^j{P_k(x)P_k(y)},\eqno (20)
$$
$$
Wr(P_{i-1}(q),P_i(q))+Wr(P_j(q),P_{j+1}(q))=(1-x^{-2})\sum_{k=i}^j{P_k^2(x)},\eqno (21)   
$$
Такие же формулы справедливы, если заменить все $P_i(q)$ на $Z_i(q)$ и $q^{-1}$ на $q^{-1}(1-q^a)(1-q^b)$. Имеют место матричные варианты всех этих тождеств, вытекающие из формул~(11) и (12).

  4) Для $0\le i\le n$
$$
qP_i(q)=\sum_{0\le k}{d_{i0}^kz^{-k-1}}.
$$
\end{prop}
 
  Из формулы~(17) вытекает, что цепные дроби из Примеров~1 и 4 равны $qP_0(q)$ для соответствующей аффинной диаграммы Дынкина. Те же формулы можно получить, если использовать в индуктивной процедуре разложения в цепную дробь не формулу дополнения по Шуру, а уравнения любой из систем~(14)-(16). При этом для аффинной диаграммы Дынкина $\tilde A_n$  получится бесконечная ветвящаяся цепная дробь. Чтобы получить цепную дробь из Примера~4, нужно остановить индуктивную процедуру на $[n/2]$-м шаге. Начиная с любого промежуточного шага в индуктивной процедуре разложения в цепную дробь для $qP_0(q)=P_0(q)/P_{-1}(q)$, мы получим разложения для всех отношений $P_{i{\to }}(q)/P_{i}(q)$. Другими словами, цепные дроби для отношений $P_{i{\to }}(q)/P_{i}(q)$ являются частями цепной дроби для отношения $P_0(q)/P_{-1}(q)$. Аналогичное утверждение справедливо для тождеств Кристоффеля-Дарбу.

  Пример 5. Испольуя формулу~(17) и Пример~3, приведём для аффинных диаграмм Дынкина $\tilde D_n$, $\tilde E_6$, $\tilde E_7$, $\tilde E_8$ списки отношений рядов Пуанкаре $P_{j_{\to }}(q)/P_j(q)$, $j\ge 0$. Некоторые из этих отношений в п.~4 будут выражены через инварианты Милнора стринг-зацеплений.
$$
\tilde D_n:
$$
$$
\frac{D_n(q)}{\tilde D_n(q)}, \frac{A_1(q)D_{n-2}(q)}{D_n(q)},\frac{D_{n-3}(q)}{D_{n-2}(q)},\dots,\frac{D_4(q)}{D_5(q)},\frac{A_3(q)}{D_4(q)},\frac{A_1^2(q)}{A_3(q)},\frac{1}{A_1(q)},\frac{1}{A_1(q)},\frac{1}{A_1(q)}
$$
$$
\tilde E_6:
$$
$$
\frac{E_6(q)}{\tilde E_6(q)}, \frac{A_5(q)}{E_6(q)},\frac{A_2^2(q)}{A_5(q)},\frac{A_1(q)}{A_2(q)},\frac{A_1(q)}{A_2(q)},\frac{1}{A_1(q)},\frac{1}{A_1(q)}
$$
$$
\tilde E_7:
$$
$$
\frac{E_7(q)}{\tilde E_7(q)},\frac{D_6(q)}{E_7(q)},\frac{A_5(q)}{D_6(q)},\frac{A_1(q)A_3(q)}{A_5(q)},\frac{A_2(q)}{A_3(q)},\frac{A_1(q)}{A_2(q)},\frac{1}{A_1(q)},\frac{1}{A_1(q)}
$$
$$
\tilde E_8:
$$
$$
\frac{E_8(q)}{\tilde E_8(q)},\frac{E_7(q)}{E_8(q)},\frac{E_6(q)}{E_7(q)},\frac{D_5(q)}{E_6(q)},\frac{A_4(q)}{D_5(q)},\frac{A_1(q)A_2(q)}{A_4(q)},\frac{A_1(q)}{A_2(q)},\frac{1}{A_1(q)},\frac{1}{A_1(q)}
$$

  \section {Ряды Пуанкаре и перестройки зацеплений}

  Стринг-зацеплением $S$ называется упорядоченный  и не имеющий особенностей набор из $n$ кривых $s_i$ в прямом произведении диска $D^2$ на отрезок $I=[-1,1]$, причём, кривые соединяют точки $(p_1,-1),\dots, (p_n,-1)$ с точками $(p_1,1),\dots, (p_n,1)$, соответственно. Соединяя для каждого $i$ точку $(p_i,-1)$ с точкой $(p_i,1)$ кривой, идущей по границе цилиндра $D^2\times I$, мы получим зацепление из $n$ компонент ${\tilde s}_i$. Назовём это зацепление вертикальным замыканием стринг-зацепления $S$. Горизонтальным замыканием двухкомпонентного стринг-зацепления назовём узел, полученный соединением точек $(p_1,j)$ и $(p_2,j)$, $j=-1,1$. Сформулируем теорему Левина \cite{19} для двухкомпонентного стринг-зацепления -- нам понадобится только этот частный случай. Обозначим через $H$ и $V$ горизонтальное и вертикальное замыкания двухкомпонентного стринг-зацепления $S$. Тогда для  соответствующих многочленов Александера-Конвея $A_H(q)$ и $A_V(q)$ имеет место формула
$$
  \frac{AV(q)}{AH(q)}=(u+1)^{1/2}\sum_{0\le k}{(\sum_{i_1,\dots, i_k}\mu_{i_1,\dots, i_k,1,1}(S))u^{k+1}},
$$
где $\mu_{i_1,\dots, i_m)}(S)$ -- инварианты Милнора стринг-зацепления $S$, $q-q^{-1}=u/(u+1)^{1/2}$. Правая часть этого равенства выражается также через индекс зацепления на бесконечном циклическом накрытии узла \cite{30}.

  Приведём краткое определение инвариантов Милнора стринг-зацепления (подробности можно найти в \cite{19}). Обозначим через $\pi$ фундаментальную группу пространства $(D^2\times I)\setminus S$. Реализуем свободную группу $F_n$ с образующими $x_1,\dots, x_n$ как фундаментальную группу пространства $(D^2\times 0)\setminus {\{p_1,\dots, p_n\}}$, причём таким образом, что для всех $i$ индекс зацепления петли $x_i$ и нити $s_i$ равен $1$. В.~Магнус определил вложение нильпотентного пополнения группы $F_n$ в кольцо рядов от $n$ некоммутирующих переменных $Z[[u_1,\dots, u_n]]$ (образующей $x_i$ отвечает ряд $1+u_i$). Для каждого $k$ изоморфны факторгруппы групп $F_n$ и $\pi$ по $k$-м членам их нижних центральных рядов. Поэтому вложение Магнуса определяет вложение нильпотентного пополнения группы $\pi$ в $Z[[u_1,\dots, u_n]]$. Для каждого $i$ можно определить элемент фундаментальной группы, представляющая петля которого идёт параллельно компоненте ${\tilde s}_i$ и её индекс зацепления с суммой всех ${\tilde s}_j$ равен $0$. Коэффициенты образов этих элементов в $Z[[u_1,\dots, u_n]]$
$$
1+\sum_{i_1,\dots, i_r}{\mu_{i_1,\dots, i_r, i}(S)u_{i_1}\dots u_{i_r}},1\le i\le n,
$$
и называются инвариантами Милнора стринг-зацепления.

  Предположим, что для некоторого диска $D^2_0\subset D^2$ пересечение зацепления $M\subset {D^2\times I}$ с цилиндром $D^2_0\times I$ совпадает с двумя параллельными хордами в $D^2_0\times 0$. Обозначим через $N$ зацепление, полученное из $M$ следующим образом: разрежем хорды перпендикулярной плоскостью (рис.~$1b$) и склеим компоненты с одной и с другой сторон разреза. Также мы можем приклеить компоненты к верхнему и нижнему основаниям цилиндра $D^2_0$ как на рис.~$1a$ или как на рис.~$1c$. В следующих случаях получатся стринг-зацепления в смысле нашего определения (обозначим их через $M_1$ и $M_2$): 

1) $M$ является двухкомпонентным зацеплением, хорды в диске $D^2_0\times 0$ принадлежат разным компонентам и приклеивание происходит как на рис.~$1a$ (очевидно, что горизонтальное замыкание стринг-зацепления $M_1$ гомотопно $N$); 

2) $M$ является узлом (однокомпонентным зацеплением) и приклеивание происходит как на рис.~$1c$ (очевидно, что горизонтальное замыкание стринг-зацепления $M_2$ гомотопно $M$).

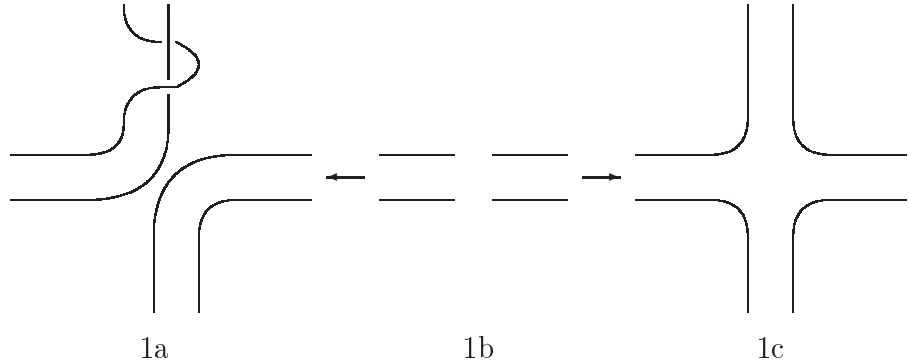
\begin{picture}(120,50)
\put(0,22){\line(1,0){10}}
\put(0,28){\line(1,0){10}}
\qbezier(10,28)(15,28)(15,32)
\qbezier(10,22)(21,22)(21,32)
\put(21,32){\line(0,1){4}}
\qbezier(15,32)(15,37)(20,37)
\put(20,37){\line(1,0){2}}
\qbezier(22,37)(28,40)(22,43)
\qbezier(20,43)(15,43)(15,48)
\put(21,38){\line(0,1){10}}
\put(17,1){1a}
\put(30,22){\line(1,0){10}}
\put(30,28){\line(1,0){10}}
\qbezier(30,28)(19,28)(19,17)
\qbezier(30,22)(25,22)(25,17)
\put(19,17){\line(0,-1){10}}
\put(25,17){\line(0,-1){10}}
\put(47,25){\vector(-1,0){5}}
\put(49,22){\line(1,0){10}}
\put(49,28){\line(1,0){10}}
\put(64,22){\line(1,0){10}}
\put(64,28){\line(1,0){10}}
\put(76,25){\vector(1,0){5}}
\put(60,1){1b}
\put(83,22){\line(1,0){10}}
\put(83,28){\line(1,0){10}}
\put(109,22){\line(1,0){10}}
\put(109,28){\line(1,0){10}}
\qbezier(93,22)(98,22)(98,17)
\qbezier(93,28)(98,28)(98,33)
\qbezier(109,22)(104,22)(104,17)
\qbezier(109,28)(104,28)(104,33)
\put(98,17){\line(0,-1){10}}
\put(98,33){\line(0,1){15}}
\put(104,17){\line(0,-1){10}}
\put(104,33){\line(0,1){15}}
\put(99,1){1c}
\end{picture}

 Обозначим через $W$ образ иммерсии набора отрезков и окружностей в диск $D^2$, причём, в качестве особенностей допускаются только двойные трансверсальные пересечения, в том числе и с границей диска. Н.~А'Кампо сопоставил кривой $W$ зацепление. Конструкция Н.~А'Кампо понадобится нам в варианте \cite{11}, точнее, нам достаточно будет знать, что соответствующее кривой зацепление реализуется в окрестности диска $D^2\times 0$ в цилиндре $D^2\times I$ как граница идущей вдоль кривой ленты с некоторыми перекрутками вдали от двойных точек и указанием в окрестности двойной точки какая компонента границы ленты проходит выше и какая ниже. В окрестности точки пересечения кривой с границей круга концы границы ленты склеиваются. Kаждый образ отрезка даёт одну компоненту в зацеплении, а каждый образ окружности -- две компоненты. В окрестности неособой точки $w$ кривой $W$ мы можем продеформировать ленту в часть плоскости диска $D^2\times 0$ и получить ситуацию, описанную в предыдущем абзаце. Более того, если точку $w$ можно соединить с границей диска $D^2\times 0$ кривой $c_w$, пересекающей кривую $W$ только в точке $w$ (назовём такую точку внешней), то в случае~1) вертикальное замыкание стринг-зацепления $M_1$ гомотопно $M$, а в случае~2) вертикальное замыкание стринг-зацепления $M_2$ гомотопно $N$. Если мы разрезаем кривую $W$ во внешней точке $w$, то концы разреза можно соединить с границей диска $D^2\times 0$ кривыми, которые параллельны кривой $c_w$.  Новой кривой, полученной таким способом из кривой $W$, конструкция Н.~А'Кампо сопоставит зацепление, которое гомотопно зацеплению $N$. На рис.~$1$ для внешней точки разреза предполагается, что большая часть кривой находится перед плоскостью страницы (в перпендикулярном ей диске $D^2$), а кривая $c_w$ находится за плоскостью страницы и там же происходит вертикальное замыкание. Из теоремы Левина вытекает

\begin{prop}\label{prop3}
Предположим, что точка разреза $w\in W$ является внешней. Если кривая $W$ удовлетворяет условиям случая~1), то
$$
   \frac{A_M(q)}{A_N(q)}=(u+1)^{1/2}\sum_{0\le k}{(\sum_{i_1,\dots, i_k}\mu_{i_1,\dots, i_k,1,1}(M_1))u^{k+1}}.
$$
Если кривая $W$ удовлетворяет условиям случая~2), то
$$
   \frac{A_N(q)}{A_M(q)}=(u+1)^{1/2}\sum_{0\le k}{(\sum_{i_1,\dots, i_k}\mu_{i_1,\dots, i_k,1,1}(M_2))u^{k+1}}.
$$
В частности, следующие отношения многочленов Кокстера выражаются через инварианты Милнора соответствующих стринг-зацеплений (на диаграммах обозначена точка разреза $w$):

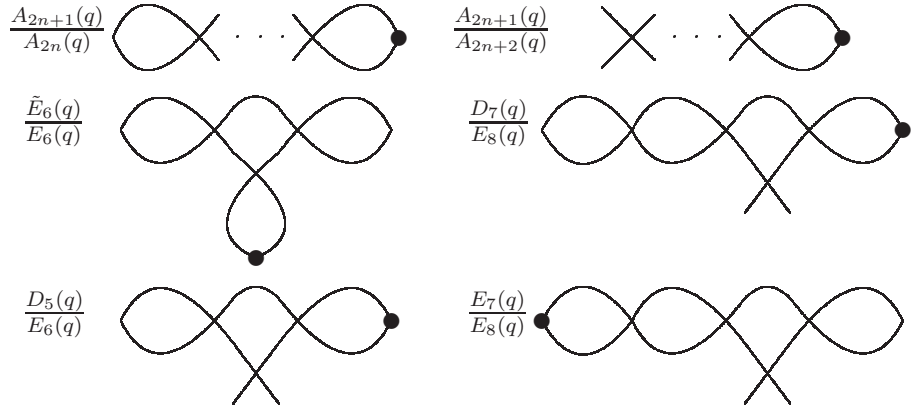
\begin{picture}(120,21)
\put(0,10){$\frac{A_{2n+1}(q)}{A_{2n}(q)}$}
\qbezier(28,7)(18,20)(14,10)
\qbezier(28,13)(18,0)(14,10)
\multiput(30,10)(3,0){3}{.}
\qbezier(38,7)(48,20)(52,10)
\qbezier(38,13)(48,0)(52,10)
\put(52,10){\circle*{2}}
\put(59,10){$\frac{A_{2n+1}(q)}{A_{2n+2}(q)}$}
\qbezier(79,6)(83,10)(86,13)
\qbezier(79,14)(83,10)(86,7)
\multiput(88,10)(3,0){3}{.}
\qbezier(96,7)(106,20)(111,10)
\qbezier(96,13)(106,0)(111,10)
\put(111,10){\circle*{2}}
\end{picture}

\begin{picture}(120,22)
\put(2,20){$\frac{\tilde E_6(q)}{E_6(q)}$}
\qbezier(30,17)(20,30)(15,20)
\qbezier(30,23)(20,10)(15,20)
\qbezier(30,23)(33,26)(36,23)
\qbezier(36,17)(46,30)(51,20)
\qbezier(36,23)(46,10)(51,20)
\qbezier(30,17)(42,7)(33,3)
\qbezier(36,17)(24,7)(33,3)
\put(33,3){\circle*{2}}
\put(61,20){$\frac{D_7(q)}{E_8(q)}$}
\qbezier(98,17)(88,30)(83,20)
\qbezier(98,23)(88,10)(83,20)
\qbezier(71,20)(76,29)(83,20)
\qbezier(71,20)(76,11)(83,20)
\qbezier(98,23)(101,26)(104,23)
\qbezier(104,17)(114,30)(119,20)
\qbezier(104,23)(114,10)(119,20)
\qbezier(98,17)(104,9)(104,9)
\qbezier(104,17)(98,9)(98,9)
\put(119,20){\circle*{2}}
\end{picture}

\begin{picture}(120,16)
\put(2,11){$\frac{D_5(q)}{E_6(q)}$}
\qbezier(30,8)(20,21)(15,11)
\qbezier(30,14)(20,1)(15,11)
\qbezier(30,14)(33,17)(36,14)
\qbezier(36,8)(46,21)(51,11)
\qbezier(36,14)(46,1)(51,11)
\qbezier(30,8)(36,0)(36,0)
\qbezier(36,8)(30,0)(30,0)
\put(51,11){\circle*{2}}
\put(61,11){$\frac{E_7(q)}{E_8(q)}$}
\qbezier(98,8)(88,21)(83,11)
\qbezier(98,14)(88,1)(83,11)
\qbezier(71,11)(76,20)(83,11)
\qbezier(71,11)(76,2)(83,11)
\qbezier(98,14)(101,17)(104,14)
\qbezier(104,8)(114,21)(119,11)
\qbezier(104,14)(114,1)(119,11)
\qbezier(98,8)(104,0)(104,0)
\qbezier(104,8)(98,0)(98,0)
\put(71,11){\circle*{2}}
\end{picture}

Пересечение этого списка со списком Примера~5 даёт отношения рядов Пуанкаре (или обратные к ним), которые выражаются через инварианты Милнора стринг-зацеплений.
\end{prop}

  Два возможных сглаживания кривой $W$ в точке двойного самопересечения можно реализовать перестройкой зацепления: разрезаем соответствующие этой точке два участка ленты и затем склеивая нужным образом образовавшиеся концы границы ленты. Допуская кривые $W$, в которых образы границ отрезков могут лежать внутри диска $D^2\times 0$ \cite{8}, мы можем также реализовать перестройкой зацепления разрезание кривой $W$ в любой точке.

  Аналогично можно применить теорему Левина к конструкции Е.~Хиронаки \cite{12}, сопоставляющей зацепление упорядоченному набору хорд в диске. Зацепление является границей поверхности, полученной приклеиваем перекрученной ленты к концам каждой хорды (эта операция известна как плюмбинг Хопфа и является частным случаем более общей операции -- суммы Мурасуги). Порядок на хордах определяет, какая из двух лент расположена выше. Мы можем разрезать любую ленту и по обе стороны разреза склеить концы границы (аналог зацепления $N$). Для набора хорд эта операция соответствует удалению хорды, отвечающей разрезанной ленте. Сглаживанию набора хорд в двойной точке также соответствует разрезание и переклейка двух соответствующих лент. Диаграмма Дынкина набора хорд являющаяся графом пересечений набора: вершины диаграммы отвечают хордам, рёбра -- точкам пересечения хорд. Любое дерево, в частности, любая из евклидовых и аффинных диаграмм Дынкина из Предложения~3, реализуемо набором хорд \cite{12}.

  Пример 6. Рассмотрим в качестве кривой $W$ окружность без двойных точек. Конструкция А'Кампо сопоставпяет ей зацепление Хопфа (два тривиальных узла, зацепленных простейшим способом), многочлен Александера которого равен $q^{-1}-q$. Разрезанной окружности соответствует тривиальный узел, который совпадает с горизонтальным замыканием стринг-зацепления, полученного из зацепления Хопфа по схеме рис. $1a$. Инварианты Милнора этого стринг-зацепления равны (пример 8 в \cite{10}) $\mu_{i_1,\dots, i_k, 1,1)}=(-1)^{k+1}$, если $i_1=\dots =i_k=1$, и равны $0$ в других случаях. Приходим к очевидному равенству
$$
  q^{-1}-q=-u(u+1)^{1/2}(1-u+u^2-u^3+\dots ).
$$ 
С точки зрения теории особенностей перестраивается зацепление особенности $A_1:x^2+y^2=0$ в (тривиальный) узел неособой точки. Той же особенности в других координатах $A_1:x^2-y^2=0$ соответствует кривая $W$, состоящая из двух пересекающихся отрезков. Сглаживание этой кривой в двойной точке реализует перестройку зацепления Хопфа в два тривиальных узла. Можно также исходить из конструкции Е.~Хиронаки, сопоставляющей набору из одной хорды зацепление Хопфа. Удаление хорды реализует перестройку зацепления Хопфа в тривиальный узел.

  Пример 7. Особенности $A_2:x^2-y^3=0$ отвечает кривая $W$, являющаяся образом отрезка и имеющая одну двойную точку. Разрезание этой кривой в неособой точке на петле реализуют перестройку узла этой особенности (трилистника) в узел особенности $A_1$ (зацепление Хопфа). Два сглаживания кривой $W$ в двойной точке реализуют перестройку трилистника в тривиальный узел или в объединение тривиального узла и зацепления Хопфа. Конструкция Е.~Хиронаки сопоставляет трилистник паре пересекающихся хорд. Удаление хорды реализует перестройку трилистника в зацепление Хопфа.

  Пример 8. По кривой $W$ определяется диаграмма Дынкина \cite{2}, для которой матрица $qS+q^{-1}S^t$ имеет вид (мы обозначаем единичные матрицы разных размеров одинаково)
$$
\begin{pmatrix}
zE&-qA&qC\\
-q^{-1}A^t&zE&-qB\\
q^{-1}C^t&-q^{-1}B^t&zE
\end{pmatrix},    
$$
где $A$, $B$, $C$ -- целочисленные матрицы с неотрицательными элементами. В случае особенностей плоская кривая является линией уровня вещественной функции от двух переменных, вершины диаграммы Дынкина соответствуют максимумам, минимумам и сёдлам этой функции и матрицы $A$, $B$, $C$ определяются по взаимному расположению максимумов, минимумов и сёдел. Веса рёбер диаграммы Дынкина удовлетворяют соотношениям Гусейн-Заде, извлекаемым из геометрии плоской кривой. Эти соотношения можно собрать в матричное равенство $AB=2C$. Разрезание кривой $W$ во внешней точке отвечает выбрасыванию вершины из диаграммы Дынкина и вычёркиванию соответствующих строки и столбца в матрице $qS+q^{-1}S^t$. Поэтому формула дополнения по Шуру даёт формулу факторизации -- связывает многочлены Кокстера исходной кривой и кривой с разрезом. Аналогичный факт имеет место для нескольких последовательных разрезов (точка может не быть внешней для исходной кривой, но стать такой после нескольких разрезов). Приведём пример такой факторизации. Существуют последовательности разрезов, приводящие к следующим матрицам 
$$
\begin{pmatrix}
zE&-qA\\
-q^{-1}A^t&zE
\end{pmatrix},
\begin{pmatrix}
zE&-qB\\
-q^{-1}B^t&zE
\end{pmatrix}.  
$$
Применяя формулу дополнения по Шуру и используя соотношения Гусейн-Заде, получим (через $d$ обозначим количество трансверсальных самопересечений кривой)
$$
\frac{G(q)}{|zE-z^{-1}AA^t||zE-z^{-1}B^tB|}=z^d|E-(4z^{-2}-1)C^t(zE-z^{-1}AA^t)^{-1}C(zE-z^{-1}B^tB)^{-1}|.
$$ 
Интересно было бы выразить правую часть этого равенства через инварианты Милнора \cite{19}, $q$-индекс зацепления \cite{30} и представление Бурау для стринг-зацеплений \cite{14}.
 
  По разрезанному зацеплению $M$ (рис.~$1b$, рис.~$2b$) можно построить стринг-зацепления другим способом:

\begin{picture}(120,47)
\put(0,22){\line(1,0){10}}
\put(0,28){\line(1,0){10}}
\put(26,22){\line(1,0){10}}
\put(26,28){\line(1,0){10}}
\qbezier(10,22)(13,22)(13,24)
\qbezier(14,27)(15,27)(15,33)
\qbezier(18,37)(15,36)(15,33)
\qbezier(18,37)(21,38)(21,41)
\qbezier(10,28)(15,28)(15,17)
\qbezier(26,22)(21,22)(21,17)
\qbezier(26,28)(21,28)(21,33)
\qbezier(19,36)(21,35)(21,33)
\qbezier(17,38)(15,39)(15,41)
\put(15,17){\line(0,-1){10}}
\put(21,17){\line(0,-1){10}}
\put(15,41){\line(0,1){4}}
\put(21,41){\line(0,1){4}}
\put(16,1){2a}
\put(45,25){\vector(-1,0){5}}
\put(45,22){\line(1,0){10}}
\put(45,28){\line(1,0){10}}
\put(62,22){\line(1,0){10}}
\put(62,28){\line(1,0){10}}
\put(74,25){\vector(1,0){5}}
\put(58,1){2b}
\put(83,22){\line(1,0){10}}
\put(83,28){\line(1,0){10}}
\put(109,22){\line(1,0){10}}
\put(109,28){\line(1,0){10}}
\qbezier(93,22)(96,22)(96,24)
\qbezier(97,27)(98,27)(98,33)
\qbezier(93,28)(98,28)(98,17)
\qbezier(109,22)(104,22)(104,17)
\qbezier(109,28)(104,28)(104,33)
\put(98,17){\line(0,-1){10}}
\put(98,33){\line(0,1){12}}
\put(104,17){\line(0,-1){10}}
\put(104,33){\line(0,1){12}}
\put(99,1){2c}
\end{picture}

  Стринг-зацепление на рис.~$2a$ получено из стринг-зацепления на рис.~$2c$ умножением на стандартную образующую группы кос. Если зацепление $M$ построено по плоской кривой $W$ и разрез соответствует внешней точке кривой, то вертикальное замыкание стринг-зацепления на рис.~$2a$ ($2c$) гомотопно $M$ ($N$). Поэтому формула~(1) выражает отношение $A_N(q)/A_M(q)$ через представление Бурау. Для разрезанной ленты в конструкции Е.~Хиронаки всё аналогично.

\bigskip

\end {document}